\documentclass[final]{siamltex1213}

\newcommand{\doingrevtex}[1]{}

\usepackage{amsmath}
\usepackage{graphicx}
\usepackage{bmpsize}
\usepackage{dsfont}
\usepackage{algorithm}
\usepackage[noend]{algpseudocode}
\usepackage{amsfonts}

\newcommand*\Let[2]{\State #1 $\gets$ #2}

\algrenewcommand\algorithmicrequire{\textbf{Precondition:}}
\algrenewcommand\algorithmicensure{\textbf{Postcondition:}}

\usepackage{breakcites}
\usepackage{microtype}
\usepackage{epstopdf}

\newcommand{\mat}[1]{\boldsymbol{#1}}

\doingrevtex{

\theoremstyle{plain}

\theoremstyle{remark}

\theoremstyle{plain}

\providecommand{\lemmaname}{Lemma}
\providecommand{\propositionname}{Proposition}
\providecommand{\remarkname}{Remark}
\providecommand{\theoremname}{Theorem}
}

\begin{document}
\title{A N-Body Solver for Free Mesh Interpolation
\thanks{This article was released under LA-UR-15-28319.  The Los Alamos National
Laboratory is operated by Los Alamos National Security, LLC for the NNSA of the
USDoE under Contract No.  DE-AC52- 06NA25396.}}
\author{Matt Challacombe\footnotemark[1] \footnotemark[2]}

\maketitle
\renewcommand{\thefootnote}{\fnsymbol{footnote}}
\footnotetext[1]{Los Alamos National Laboratory}
\footnotetext[2]{matt.challacombe@freeon.org}
\renewcommand{\thefootnote}{\arabic{footnote}}
\pagestyle{myheadings}
\thispagestyle{plain}
\markboth{Challacombe}{A N-Body Solver for Free Mesh Interpolation}

\begin{abstract}
Factorization of the Gaussian RBF kernel is developed for 
free-mesh interpolation in the flat, polynomial limit corresponding to 
Taylor expansion and the Vandermonde basis of geometric moments.
With this spectral approximation, a top-down octree-scoping of 
an interpolant is found by recursively decomposing the residual,
similar to the work of Driscoll and Heryudono (2007), except that in the current approach
the grid is decoupled from the low rank approximation, allowing partial separation of 
sampling errors (the mesh) from representation errors (the polynomial order).
Then, it is possible to demonstrate roughly 5 orders of magnitude improvement in 
free-mesh interpolation errors for the three-dimensional Franke function, relative to previous benchmarks. 
As in related work on $N$-body methods for factorization by 
square root iteration (Challacombe 2015), some emphasis is placed on resolution of the identity.
\end{abstract}

\maketitle
\section{Introduction}

A reproducing kernel Hilbert space (RKHS) ${\cal H}_d$  is associated with $d$-dimensional kernels $K_d$,  
mapping $\mathbb{R}^d \times \mathbb{R}^d \rightarrow \mathbb{R}$ via an inner product $\left<\,\cdot\, , \,\cdot\, \right>_{\cal H}$
\cite{steinwart2006explicit,hofmann2008kernel,schaback2009nonstandard,fasshauer2011positive}.
Kernels of the RKHS are reproducing: 
\begin{equation*}
{K}(\mat{x}',\mat{x}) =  \left< K(\,\cdot\,, \mat{x}' ) , K(\,\cdot\,, \mat{x}) \right> \, ,
\end{equation*}
with the ability to reconstruct all functions in the native space:
\begin{equation*}
f({\mat{x}}) = \left< f(\,\cdot\,) , K(\,\cdot\,,\mat{ x}) \right> \, .
\end{equation*}
The natural bases for the RKHS are the kernel eigenfunctions $\varphi_i$:
\begin{eqnarray*}
{K}(\mat{x}',\mat{x}) &=& \sum_i \lambda_i \varphi_i \left(\mat{x}'\right) \varphi_i \left(\mat{x}\right)\\
         &=&  \left< \mat{\sqrt{\lambda}\varphi}' , 
\mat{\sqrt{\lambda} \varphi} \right> \,.
\end{eqnarray*}
In many applications, the RKHS provides an extremely powerful framework for interpolation, reconstruction, convolutions \& so on. 
Startling advances have been achieved for higher dimensional problems associated with inference and probability,
through new formalisms that are asymptotically well behaved with respect to dimensionality 
\cite{Belkin2002,griebel2005sparse,burges2009dimension,weinberger2004learning,fasshauer2012average,fasshauer2012dimension}.  

In practical applications, the pointwise function $\mat{f}$ and the corresponding Kernel or Gramian matrix $\mat{K}$ 
determine a function through the linear system $\mat{K}\cdot \mat{\omega} = \mat{f}$: 
 \begin{equation*}
{f}(\mat{x}) = \sum_{i} \omega_i \, K \left( \mat{x} , \mat{x}_i \right) \, .
\end{equation*}
Kernels like the Gaussian Radial Basis Function kernel,
\begin{equation*}
\phi \left(\lVert \mat{x}'-\mat{x} \rVert \right) = e^{-\varepsilon^2 \left( \mat{x}'-\mat{x}\right)^2 }\, ,
\end{equation*}
have many excellent properties, including translational invariance, a separable $d$-dimensional product form,  
analytic solutions to the kernel eigenfunctions \cite{fasshauer2011positive}
and super-spectral (exponential) convergence \cite{fornberg2005accuracy,platte2011fast,platte2011impossibility}. 
Another remarkable feature of the Gaussian RBF kernel is that towards the flat limit of the shape parameter, $\varepsilon \rightarrow 0$, 
interpolation tends towards approximations with early low-rank saturation but with  severely ill-conditioned
Gramian matrices, $\mat{\phi}$ \cite{maz1996approximate,driscoll2002interpolation,
                          fornberg2004some,fornberg2005accuracy,schaback2005multivariate,Fasshauer2012}.
This is the RBF-direct problem, identical to the Gaussian LCAO problem faced in quantum chemistry 
\cite{Rothlisberger2002,Jansik2007,helgaker2008molecular}. 

Little appreciated, this ill-conditioning can be partially sidesteped through only low rank 
factorization in the case of kernels with rapid spectral decay 
\cite{fasshauer2011positive,Fasshauer2012,wathen2013spectral}.  Also, recent progress in analysis of the 
Gaussian kernel shows that in the flat limit, eigenfunctions of the Gaussian RBF kernel tend to moments 
of the Hermite polynomials with rank corresponding to the number of centers 
\cite{fasshauer2011positive,Fasshauer2012,wathen2013spectral}, ${\cal L}=N_p$, equivalent to 
any other polynomial moments of the same rank, via {\em e.g.} orthogonality of the Vandermonde basis 
\cite{de2001inverse,gander2005change,o2009discrete}. The transition from exponential (Gaussian RBF) to polynomial (Vandermonde) behavior
in the flat regime gives Gaussian RBF error curves their characteristic $^{\scriptstyle \cal p}$-like dip 
\cite{fornberg2005accuracy,platte2011fast,platte2011impossibility,wathen2013spectral}.  
  
These spectral approximations may reach convergence much faster 
than a naively taken grid can establish an accurate fill distance  
\cite{DeMarchi2010,platte2011fast,zwicknagl2013interpolation}. Then low, ${\cal L}$-rank factorization  
can be much faster than the unisolvent case, ${\cal O}\left( {\cal L}^2 N_p \right)$ vs. ${\cal O}\left( N^2_p \right)$, 
allowing partial separation between complexity of the representation and of the grid.
In just 1-dimension however,  RBF expansions to numerical accuracy for challenging functions 
may require $N_p \sim 20-30$ centers \cite{Fasshauer2012}, 
corresponding to rank ${\cal L} \sim L^d$. In just three dimensions, factorization is ${\cal O}\left( L^6 N_p \right)$, already a formidable 
challenge for large (naive) grids.  Notable attacks on this problem include statistical methods for
very long and skinny matrix factorizations  \cite{halko2011finding,grasedyck2013literature,mahoney2011randomized}.

In higher dimensions and for functions that are strongly cusped, 
the number of grid points naively taken to maintain an accurate fill distance 
$h_p \sim N^{-1/d}_p$ may grow explosively \cite{wendland2004scattered,DeMarchi2010,platte2011fast,zwicknagl2013interpolation}.  
Then, factorization demands a scoping mechanism capable of culling-out  low-dimensional 
sub-space representations from higher dimensional manifolds
\cite{Belkin2002,burges2009dimension,weinberger2004learning},    
including {\em e.g.} dimension-free shape parameter sequences 
\cite{fasshauer2012dimension,fasshauer2012average} and sparse grids \cite{griebel2005sparse,Garcke2013}.  
Similarly, finding low complexity sub-structures in Gramian factorizations is
enabled by octree-scoping of the $\tt SpAMM$ bound in $n$-body approaches to square root iteration \cite{Challacombe2015}.

For potential-like RBF kernels with also spectral convergence, such as the multiquadratic, Taylor expansions converge 
rapidly only in the far-field limit. These kernels are enabled by well known methods, including the Barnes-Hut Tree-Code and the 
Fast-Multipole-Method (FMM), with octree scoping of the bi-polar expansion based on range criteria, and the 
upwards accumulation of long ranged, increasingly smooth potentials 
\cite{yokota2010petrbf,krasny2011fast,march2015askit,march2015robust}.
Different from these FMM-like approaches, the focus of this work is instead on 
top-down expansion of the residual via low-rank factorization of the Gaussian kernel in the truncated polynomial  
limit, obtaining the largest and smoothest 
components of the interpolant first, resolving increasingly fine details with decreasing magnitudes.
With large quasi-random grids, can low rank octree-scoping yield high accuracy? 

This paper is organized as follows:  In the next section, we develop polynomial iterpolation through 
Taylor expansion of the Gaussian RBF kernel, and discuss equivalence with kernel eigenfunctions in the flat limit. 
Then, we show how the low rank approximation describes the three-dimensional 
Franke function, used to benchmark other methods for scattered data interpolation 
\cite{bozzini2002testing,cavoretto2015,de2012polynomial}.  Finally, we show increasingly compact and accurate 
representation of the method with increasing expansion order and grid size, through $N_p=8^9$. 

\section{Taylor Expansions of  Gaussian Kernels}\label{indinterp}
In this section, the polynomial limit of the flat Gaussian RBF kernel is sketched in three-dimensions 
with emphasis on concepts and notation familiar to chemistry and materials science. To start, 
expansion of the Gaussian RBF kernel can be factored by polynomial order,

\begin{eqnarray*}
\phi \left(\lVert \mat{x}'-\mat{x} \rVert \right) &=& \sum_{lmn} 
\, \varepsilon^{l+m+n} \, \frac{x'^l_1 x'^m_2 x'^n_3}{l! m! n!} 
h_l\left(\varepsilon x_1 \right) h_m\left(\varepsilon x_2 \right) h_n\left(\varepsilon x_3 \right)\\
&=&\sum_{L=0} \varepsilon^L \sum^{\ell^{\,3}_{\,\scriptscriptstyle L+1}}_{lmn=\ell^{\,3}_{\,\scriptscriptstyle L}}
\, \Lambda_{lmn} \!\! \left({\mat{x}'}\right)\,  h_l\left(\varepsilon x_1 \right) h_m\left(\varepsilon x_2 \right) h_n\left(\varepsilon x_3 \right)
\end{eqnarray*}
where $h_i \left( x \right) = e^{-\varepsilon^2 x^2} H_i \left( x\right)$ are the Hermite functions, and 
\begin{equation*}
{ {{\ell}}}^{\,3}_{\, i} = i(i+1)(i+2)/6
\end{equation*}
is the rank stride of the three dimensional multi-index $lmn$.  Then, with the Taylor expansion
\begin{equation*}
h_\alpha \left( t \right) = \sum_\beta \frac{t^\beta}{\beta!} h_{\alpha+\beta} \left( 0 \right)
\end{equation*}
the Gaussian RBF kernel is 
\begin{multline*}
\phi \left(\lVert \mat{x}'-\mat{x} \rVert \right) = 
\sum_{L=0} \varepsilon^{\scriptscriptstyle L} 
\sum_{L'=0} \varepsilon^{\scriptscriptstyle L'}  
 \sum^{\ell^{\,3}_{\,\scriptscriptstyle L+1}}_{lmn=\ell^{\,3}_{\scriptscriptstyle L}} 
\; \sum^{\ell^{\,3}_{\,\scriptscriptstyle L'+1}}_{lmn'=\ell^{\,3}_{\scriptscriptstyle L'}}
\; \Lambda_{lmn} \!\! \left({\mat{x}}\right)\, 
   \Lambda_{lmn'} \! \left({\mat{x}'}\right)\, \\ \times (-1)^{l+m+n} \,  h_{l+l'}\left(0\right) h_{m+m'}\left(0\right) h_{n+n'}\left(0\right)\, .
\end{multline*}
So far, the unisolvent ${\cal L} = N_p$ case has been implied.  Now please consider a truncated expansion with
${{\cal L}\equiv  \ell^{\,3}_{\,{L_{\rm max}}} } \ll N_p$.  Then, the kernel expansion may be written:
\begin{equation} \label{fullyfactored}
\mat{\phi}_{\scriptscriptstyle N_p \times N_p} = \mat{\Lambda}_{\scriptscriptstyle N_p \times \cal{L}} 
\cdot \mat{\varepsilon}_{\cal{L} \times \cal{L} }  \cdot 
\mat{T}_{\cal{L} \times \cal{L} } 
\cdot \mat{\varepsilon}_{\cal{L} \times \cal{L} }  \cdot  
\mat{\Lambda}^T_{{\cal{L} }\times \scriptscriptstyle N_p } \; ,
\end{equation}
where 
\begin{equation*}
 \mat{\varepsilon}=diag\left[\; \varepsilon^0, \; \varepsilon^1,\varepsilon^1,\varepsilon^1, 
        \; \varepsilon^2, \dots \;,\varepsilon^{{\scriptscriptstyle L}_{\rm max}}
 \; \right]_{{\cal L}\times{\cal L}} \;,  
\end{equation*} 
\begin{equation*}
T_{lmn, lmn'} = (-1)^{l+m+n} \,  h_{l+l'}\left(0\right) h_{m+m'}\left(0\right) h_{n+n'}\left(0\right) \, ,
\end{equation*} 

\begin{equation*}
h_\alpha(0)=\left\{  
\begin{array}{cc}
0 & \alpha \qquad {\rm odd} \\
(-1)^\frac{\alpha}{2}(\alpha -1)!!  & \alpha \qquad {\rm even} 
\end{array}
\right. \; 
\end{equation*}
are the Hermite numbers, and 
 \begin{equation*}
{\Lambda}_{lmn} \left( \mat{x} \right) = \frac{x^l_1 x^m_2 x^n_3}{l! m! n!} 
\end{equation*}
are the geometric moments \cite{hu1962visual,schaback2005multivariate,o2009discrete,Flusser2009},
also known as the unabridged multipole moments in chemistry \cite{Applequist1995}.   

QR factorization of the Vandermonde matrix of polynomial moments,
\begin{equation*}
\mat{\Lambda}  = \mat{Q}_{\scriptscriptstyle N_p \times \cal{L}} 
                 \cdot \mat{R}_{\scriptscriptstyle \cal{L} \times \cal{L}}  \, ,
\end{equation*}
provides a representation-free basis, $\mat{Q}$, as shown by Fornberg and others    
\cite{gander1980algorithms,bjorck2004calculation,fornberg2007stable,Pazouki2011,Fasshauer2012},
and including the early work of Gene Golub \cite{golub1965numerical,golub1965calculating}.  

Then, 
\begin{equation*}
\mat{\phi}_{\scriptscriptstyle N_p \times N_p} = \mat{Q}_{\scriptscriptstyle N_p \times \cal{L}} 
                                    \cdot \mat{\Phi}_{\scriptscriptstyle \cal{L} \times \cal{L}} 
\cdot \mat{Q}^T_{{\cal{L} }\times \scriptscriptstyle N_p } \, ,
\end{equation*}
where the auxiliary representation of the kernel is:
\begin{equation*}
  \mat{\Phi} = \mat{R} \cdot  \mat{\varepsilon} \cdot  \mat{V}\cdot \mat{e} \cdot \mat{V}^T  \cdot  \mat{\varepsilon} \cdot   \mat{R}^T \, ,
\end{equation*}
with $\mat{V}$ and $\mat{Q}$ orthogonal,  $\mat{T}=\mat{V}\cdot \mat{e} \cdot \mat{V}^T$, and where $\mat{e}$ and $\mat{\varepsilon}$ are diagonal.
By the properties of orthogonal matrices, the inverse is simply
\begin{equation}
\mat{\phi}^{-1}  =\mat{Q} \cdot [\mat{R}^{-1}]^T \cdot  \mat{\varepsilon}^{-1} \cdot  \mat{V}\cdot \mat{e}^{-1} \cdot \mat{V}^T  
\cdot  \mat{\varepsilon}^{-1} \cdot   \mat{R}^{-1} \cdot \mat{Q}^T \; ,
\end{equation}
and relative expressions involving poorly behaved powers of the shape parameter can be handled safely in a number of ways. 
So far, this is the RBF-QR analysis of Fornberg \cite{} and others \cite{}  with 
analytically well conditioned solutions due to careful handling as $\varepsilon \rightarrow 0 $.





Sidesteping an elegant treatment of the shape-factor, consider now the problem of transforming a pointwise interpolant $p$,
\begin{equation*}
\mat{f}_p = \{  f\left(x_{i,1},x_{i,2},x_{i,3}\right) | i=1,N_p \}\, ,
\end{equation*}
sampled on a random grid, 
\begin{equation*}
\mat{x}_p = \{  [\underline{a},\underline{b},\underline{c}]_i \mid [\underline{a},\underline{b},\underline{c}] \in [-1,1]^3 , i=1,N_p \} \, 
\end{equation*}
to another randomly chosen grid, $\mat{x}_{q} \ne \mat{x}_p$.  
With $'$ denoting $q$ dependent terms, the transformation between these free-mesh pointsets $p \rightarrow q$ is:
\begin{eqnarray*}
\mat{f}'_{q} &=& \mat{\phi}'_{N_q \times N_p} \cdot \mat{\lambda}_{p}\\
                &=& \mat{\phi}'_{N_q \times N_p} \cdot \mat{\phi}^{-1}_{N_p\times N_p} \cdot \mat{f}_{p}\\
                &=& \mat{Q}'_{N_q \times {\cal L}}\cdot \mat{\Phi}'^{\,1/2}_{{\cal L} \times {\cal L}} 
                     \cdot \left[ \mat{\Phi}^{-1/2}\right]^T_{{\cal L}\times \cal L} \cdot \mat{Q}^T_{{\cal L}\times{N_p}} \cdot \mat{f}_{p}\\
                &=&  \mat{I}'_{N_q \times N_p} \cdot \mat{f}_{p} \; ,
\end{eqnarray*}
where 
\begin{equation}\label{identityprime}
\mat{I}' = \mat{Q}' \cdot [ \mat{R}']^T \cdot \mat{R}^{-1} \cdot \mat{Q}^T
\end{equation}
is the pointwise resolution of identity.

This somewhat surprising, shape-parameter-free result does not follow from a special handling 
of the shape parameter, but obtains simply from polynomial independence in the factored Taylor space, 
via orthogonality of the Vandermonde basis \cite{de2001inverse,gander2005change,o2009discrete}. 

\begin{figure}[h] 
\fbox{\includegraphics[width=5.9cm,keepaspectratio=true,
                       trim={.5cm 2.3cm 2.cm 1.cm},clip]{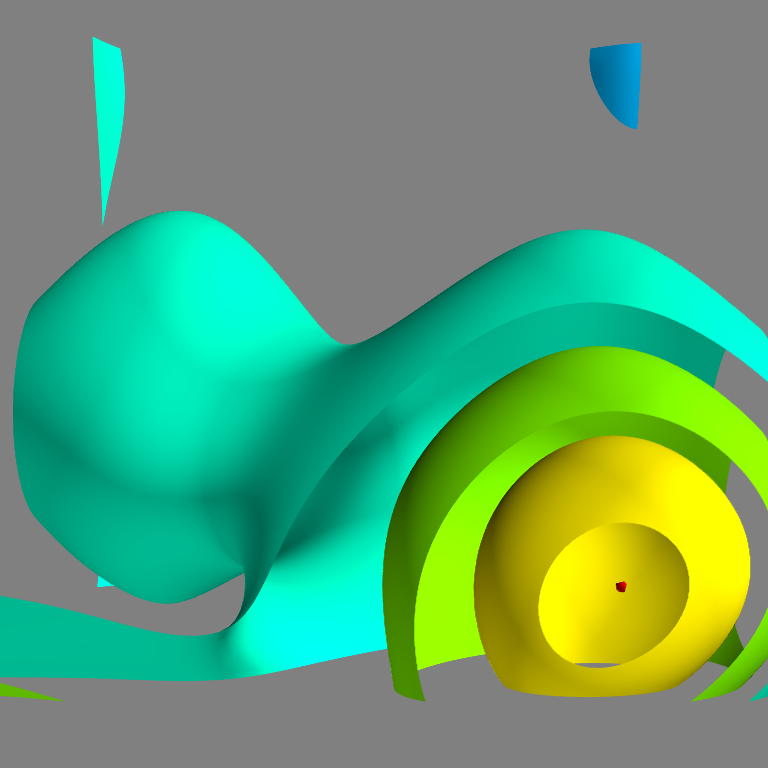}} 
\fbox{\includegraphics[width=5.9cm,keepaspectratio=true,
                       trim={.5cm 2.3cm 2.cm 1.cm},clip]{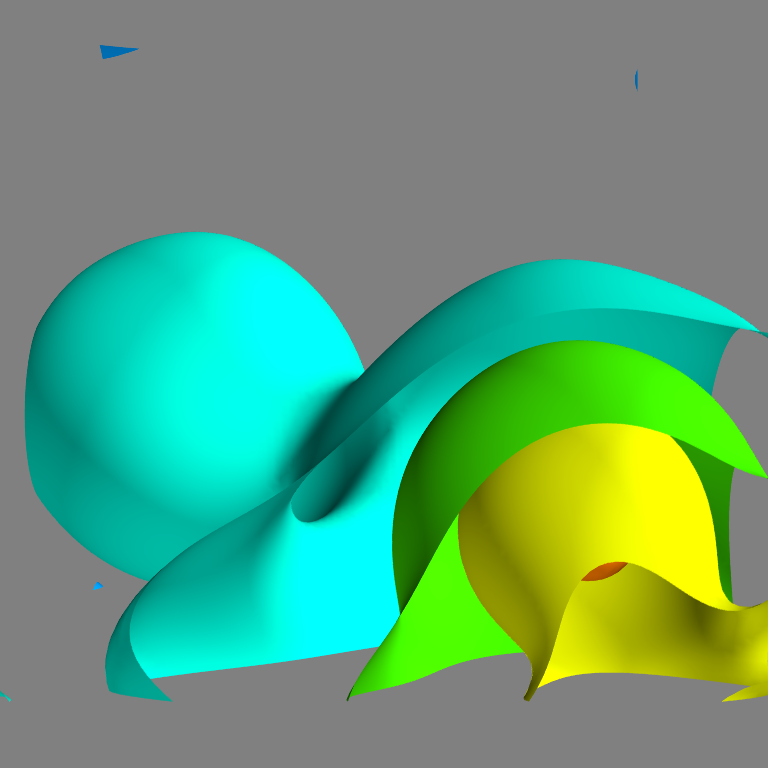}} \\
\fbox{\includegraphics[width=5.9cm,keepaspectratio=true,
                       trim={.5cm 2.3cm 2.cm 1.cm},clip]{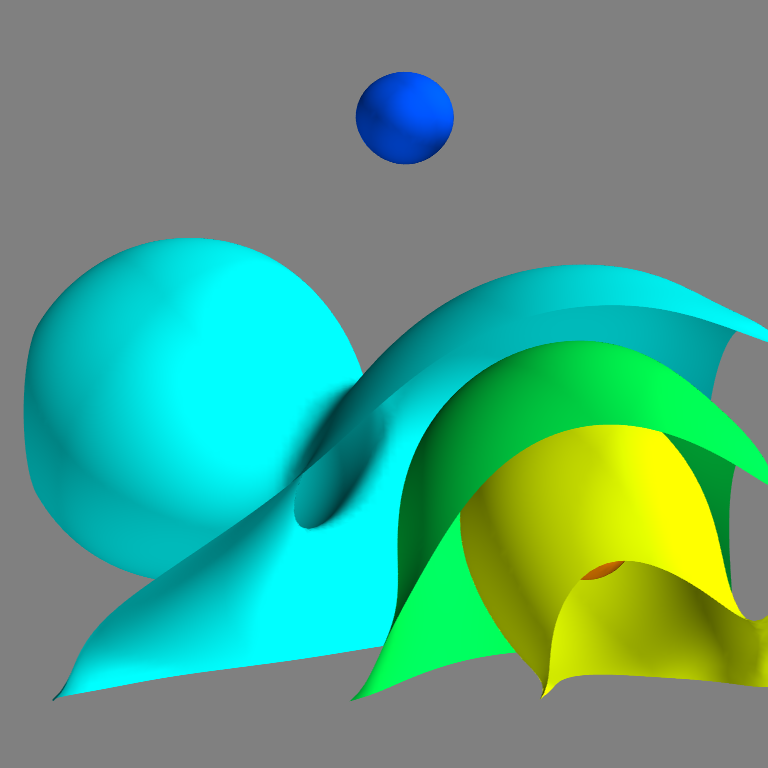}} 
\fbox{\includegraphics[width=5.9cm,keepaspectratio=true,
                       trim={.5cm 2.3cm 2.cm 1.cm},clip]{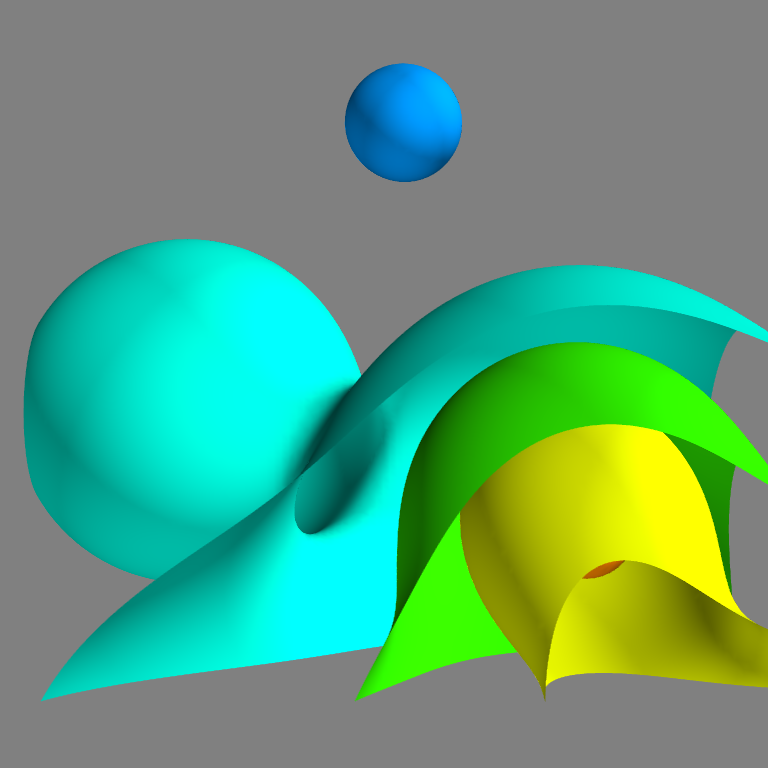}} 
\caption{
Unscoped free-mesh interpolation of the three-dimensional Franke function, $f_{\rm F3d}$, with low rank polynomial approximation
corresponding to random, $N_p=8^6$ grids in $[0,1]^3$.  Shown from left to right, top to bottom, are 
equispaced surfaces with values between $\tt -.10$ and $\tt .79$ (blue to red) for orders $L_{\rm max} = 6, 10 \; \& \; 14$,
and in the lower rightmost pannel, surfaces of the native function.}\label{3dfranke}
\end{figure}

This result is similar to the RI method developed in quantum chemistry by Yu-Cheng Zheng and Jan Alml{\"o}f \cite{zheng1993,zheng1996}, and
related to random phase Vandermonde matrices with very interesting properties \cite{tucci2011eigenvalue,debbah2011convolution}, 
to the problem of statistical and sub-random grids \cite{morokoff1994quasi,ghahramani2002bayesian}, 
and to the computation of weakly admissible meshes \cite{bos2010least} and approximate Fetke points \cite{sommariva2009computing}.  

In this contribution however, we are interested only in grid independence in the limit of small fill distances, assumed 
here to be $h_p \sim N^{-1/3}_p$.  Also, in practice there is no need to actually form $\mat{I}'$, since 
\begin{eqnarray*}
\mat{f}'_{q} &=& \mat{\Lambda}'_{N_q \times {\cal L}} \cdot  \mat{R}^{-1}_{{\cal L}\times \cal L} \cdot \mat{Q}^T_{{\cal L}\times{N_p}} \cdot \mat{f}_{p}\\
           &=& \mat{\Lambda}'_{N_q \times {\cal L}} \cdot  \mat{F}_p \, ,
\end{eqnarray*}
where $\mat{F}_p$ is the truncated Taylor representation of $\mat{f}_p$.  More generally, 
\begin{equation*}
f \left( \mat{r} \right) = \mat{\Lambda} \left(\mat{r} \right) \cdot \mat{F}_p  \,.
\end{equation*}

\section{Low Rank Approximation of the 3-D Franke Function}
In three-dimensional studies of scattered mesh interpolation for smooth functions, the  
three-dimensional Franke function
\begin{multline}\label{eqF1}
f_{\rm F3d}\left( x,y,z\right) = \frac{3}{4}\left[ e^{-\frac{1}{4} \left[(9x-2)^2+(9y-2)^2+(9z-2)^2  \right]}+e^{-\frac{1}{49}(9x+1)^2-\frac{1}{10}(9y+1)^2-\frac{1}{10}(9z+1)^2} \right]\\
+\frac{1}{2} e^{-\frac{1}{4} \left[ (9x-7)^2+(9y-3)^2+(9z-5)^2\right]}-\frac{1}{5} e^{-(9x-4)^2-(9y-7)^2-(9z-5)^2}\; ,
\end{multline}
is perhaps the most common choice for benchmarking errors \cite{franke1979critical,bozzini2002testing,cavoretto2015,de2012polynomial}.
In Fig.~\ref{3dfranke}, interpolation surfaces of $f_{\tt F3d}$ are shown for orders $L_{\rm max} = 6, 10 \; \& \; 14$, corresponding to
${\cal L}=84, 286 \; \&  \; 680$.  In this challenging example, the correct qualitative behavior 
(ability to resolve positive and negative lobes) is obtained only for $L_{\rm max}=14$.
 \begin{algorithm}
  \caption{Recursive octree decomposition of the FMT residual.}\label{alg1}
   \begin{algorithmic}[1]
    \Statex
    \Procedure{$\tt FMT\_Mesh\_to\_Tree$}{$N, L, {\cal L}, \mat{x}_{(1:N,1:3)},\mat{f}_{(1:N)}, P $}
    \Let{$\mat{x}_{\left(1:N,1:3\right)}$}{$\left(\mat{x}_{\left(1:N,1:3\right)}-P_{\mat{C}\left({1:3}\right)}\right)*P_{{\tt Scale}}$}
    \State \Call{${\tt FMT\_Mesh\_to\_Moments}$}{$N, L,{\cal L}, \mat{x}_{\left(1:N,1:3\right)}, \mat{f}_{(1:N)}, P_{F\left(1:{\cal L}\right)} $ }
    \State \Call{$\tt FMT\_Ordered\_by\_Octant$}{$ N, \mat{x}_{\left(1:N,1:3\right)}, \mat{f}_{(1:N)}, n_{\left(1:8\right)}$} 
    \Let{$j$}{$1$} 
    \For{$o \gets 1 \textrm{ to } 8$}        
          \Let{$k$}{$\sum^o_i n_i$}
          \Let{$E_{\scriptscriptstyle \rm RMS}$}{$\sqrt{\sum^k_{i=j}f^2_{i}/n_o}$} 
        \If{ $E_{\scriptscriptstyle \rm RMS} >  {\Large \tau} \;  \& \; n_o \ge {\cal L}$}
           \State \Call{{\tt FMT\_Mesh\_to\_Tree}}{$n_o, L, {\cal L}, \mat{x}_{\left(j:k,1:3\right)} ,\mat{f}_{\left(j:k\right)} , 
{\tt FMT\_New\_Octant}\left(o, P  \right)$}
        \EndIf
           \Let{$j$}{$j+k$}
       \EndFor
     \EndProcedure
   \end{algorithmic}
  \end{algorithm}

 \begin{algorithm}
  \caption{}
   \begin{algorithmic}[1]
    \Statex
    \Procedure{$\tt FMT\_Mesh\_to\_Moments$}{$N, L,{\cal L}, \mat{x}_{\left(1:N,1:3\right)},\mat{f}_{\left(1:N\right)},P_{F\left(1:{\cal L}\right)} $} 

    \For{$l \gets 0, L$}        
    \For{$m \gets 0, L-l$}        
    \For{$n \gets 0, L-l-m$}        
        \Let{$\mat{\Lambda}_{\left(1:N, \, lmn\left[l,m,n\right]\,\right)}$}{
$ {\mat{x}^{\,l}_{\left(1:N,1\right)} \mat{x}^{\,m}_{\left(1:N,2\right)} \mat{x}^{\,n}_{\left(1:N,3\right)}}/{l! m! n!} $}
    \EndFor
    \EndFor
    \EndFor
    \Let{$\{\mat{Q},\mat{R}\}$}{${\tt FMT\_DGEQRF\_DORGOR}\left(N, {\cal L},  \mat{\Lambda}_{(1:N,{\cal L})} \right)$}   
    \Let{$P_{\mat{F}(1:{\cal L})} $}{$ \mat{R}^{-1}_{(\cal L, L)} \cdot {\mat{Q}^T}_{\!({\cal L},1:N)} \cdot \mat{f}_{\left(1:N \right)}$}
    \Let{$\mat{f}_{(1:N)}$}{$\mat{f}_{(1:N)}- \mat{\Lambda}_{(1:N,1:{\cal L})} \cdot P_{\mat{F}(1:{\cal L})}$}

    \EndProcedure
  \end{algorithmic}
 \end{algorithm}

\begin{figure}[h]
\includegraphics[width=5.in]{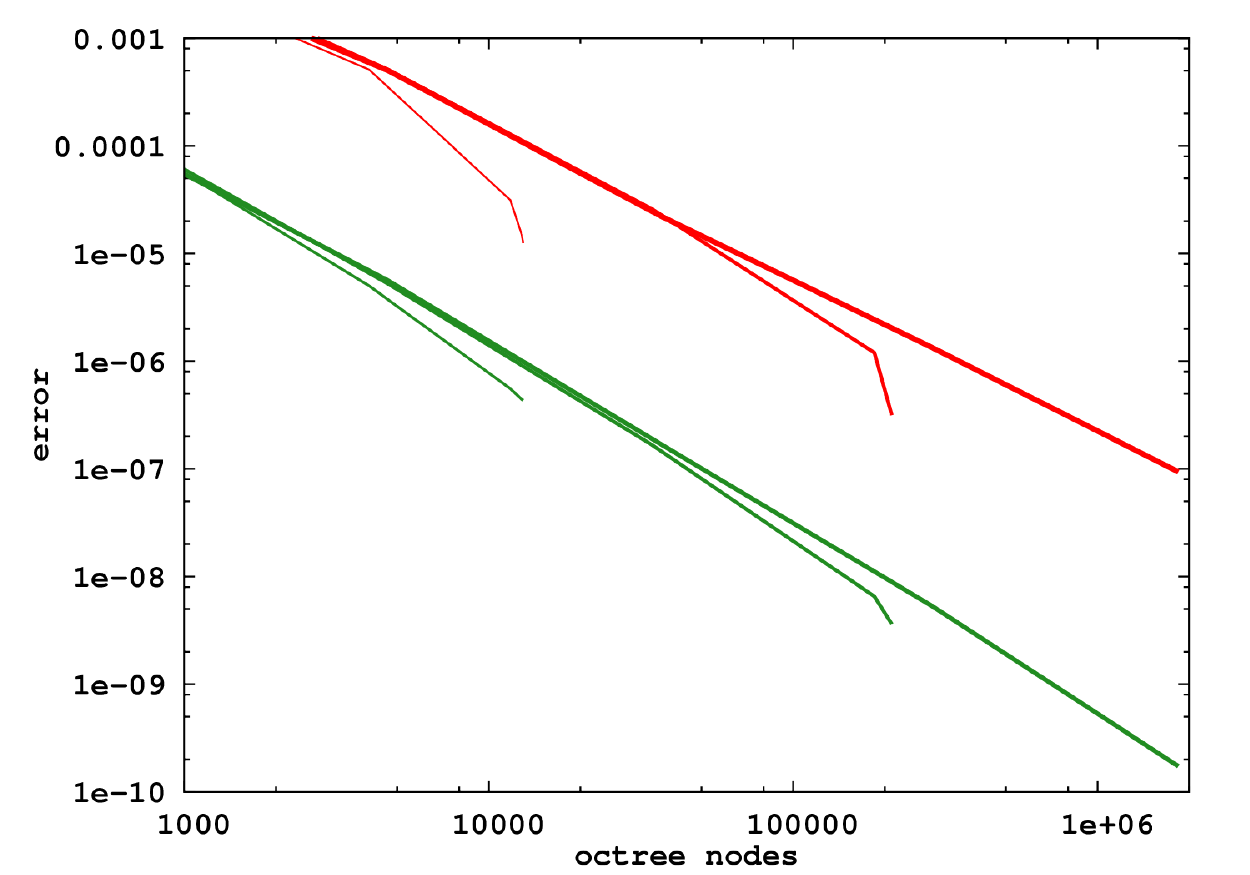}
\caption{Convergence of the $\epsilon_{\tt rms}$ and $\epsilon_{\infty}$  
errors for Fast Field Transform of $F1$, Eq.~\ref{eqF1}, on random grid $p$ ($p\in[0,1]^3$) 
to random grid $q$ ($q\in[0,1]^3$).  Shown are errors for expansion 
order ${\cal L}=6 \; \& \; 10$ with respect to the number of nodes,
corresponding to random grids of size $N_p=N_q=\{8^4, 8^5, 8^6, 8^7\}$.} \label{expt_1}
\end{figure}

\begin{figure}[h]
\includegraphics[width=5.in]{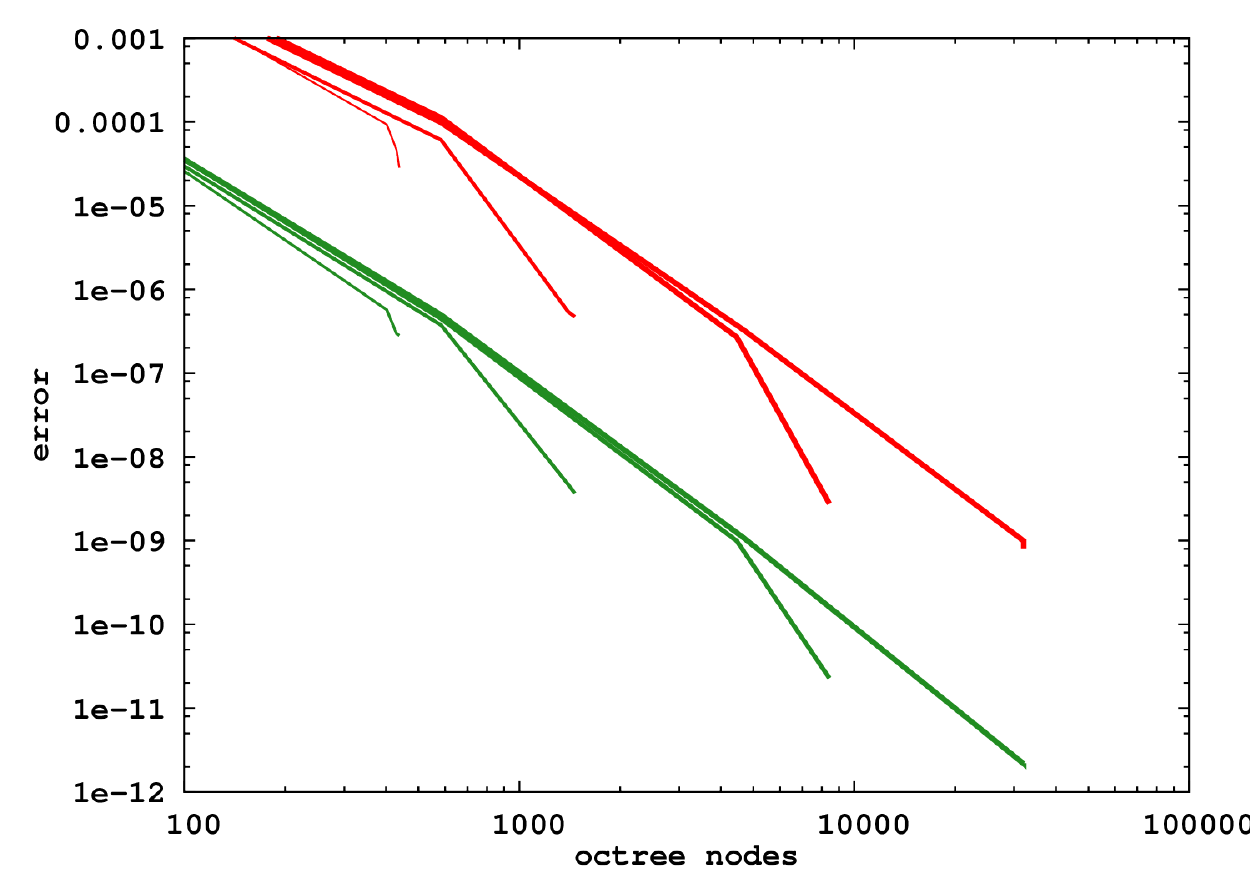}
\caption{Convergence of the $\epsilon_{\tt rms}$ and $\epsilon_{\infty}$  
errors for Fast Field Transform of $F2$, Eq.~\ref{eqF2},  on random grid $p$ ($p\in[0,1]^3$) 
to random grid $q$ ($q\in[0,1]^3$).  Shown are errors for expansion 
order ${\cal L}=6 \; \& \; 10$ with respect to the number of nodes,
corresponding to random grids of size $N_p=N_q=\{8^4, 8^5, 8^6, 8^7\}$.} \label{expt_2}
\end{figure}

\begin{figure}[h]
\includegraphics[width=5.in]{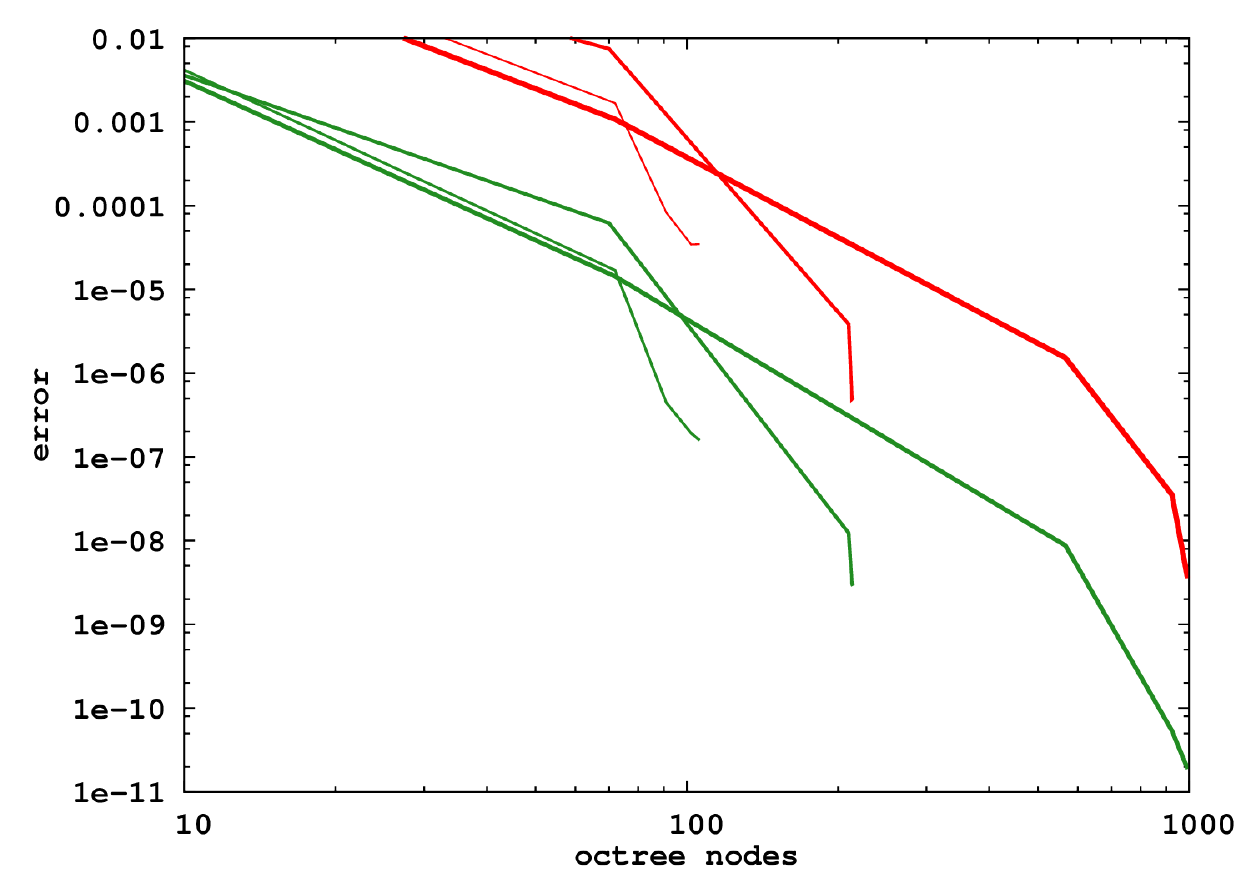}
\caption{Convergence of the $\epsilon_{\tt rms}$ and $\epsilon_{\infty}$  
errors for Fast Field Transform of $F2$, Eq.~\ref{eqF2},  on random grid $p$ ($p\in[0,1]^3$) 
to random grid $q$ ($q\in[0,1]^3$).  Shown are errors for expansion 
order ${\cal L}=6 \; \& \; 10$ with respect to the number of nodes,
corresponding to random grids of size $N_p=N_q=\{8^4, 8^5, 8^6, 8^7\}$.} \label{expt_3}
\end{figure}

\section{Free Mesh Transform}\label{fmt}
The free-mesh transform ($\tt FMT$), shown in Algorithm \ref{alg1}, is simpler than other fast methods 
with octree-scoping, with only downward recursion to accumulate local polynomial representation of the residual. 
This approach is similar to the adaptive residual sub-sampling method of Driscoll and Hryundono 
\cite{driscoll2007adaptive}, except that in the current approach, the grid is decoupled from the 
low rank approximation, allowing partial separation of sampling errors (the mesh) from representation 
errors (the Taylor order $L$).  

Starting at the top level,  the interpolant and grid, $\{\mat{f}_{N_p},\mat{x}_{(1:N_p,1:3)}\}$, together with the octree $P$,
are recursively subdivided, constructing approximations to residuals in sub-octants through coordinate translations and 
dilations that smooth and re-center the interpolant 
within $\left[-1,1\right]^3$.  At some depth, the number of points in a new octant, $n_o$, 
will approach the most effective resolution possible;  Algorithm \ref{alg1} looks ahead, terminating recursion 
for any octants obeying $n_o < {\cal L}$. 
Recursion also terminates when the local RMS error meets the convergence criteria $E_{\rm RMS}<\tau$, determined by
ability of the polynomial representation to completely annihilate local residuals.  

In Figures \ref{expt_1}-\ref{expt_3}, scoping experiments with $L_{\rm max}=4, 8, \& 12$ are carried out with thresholds
$\tau=\{10^{-6}, 10^{-8}, 10^{-10},10^{-12}\}$ and separate grids with $N_q = N_p =\{ 8^4, 8^5, 8^6, 8^7, 8^8, 8^9\}$. 
  Shown are the total RMS error;
\begin{equation*}
E_{\tt rms} = \sqrt{\frac{\sum_{i} \{ [{f}_q]_i-f(x_i,y_i,z_i) \}^2  }{N_q}}
\end{equation*}
and the max error;
\begin{equation*}
E_{\infty} = \max_i \| [f_q]_i-f(x_i,y_i,z_i) \| \, , 
\end{equation*}
with respect to the total number of octree-nodes.   These results show saturation of the representation,
set by $L_{\rm max}$, with respect to given RMS error threshold, $\tau$, and with respect to an increasing mesh, $N_q=N_p$.
Notably, with increasing order, $L_{\rm max}$, behavior of the scoping error in terms of 
octree-complexity becomes noticeably more tree like; ${\cal O}(N) \rightarrow {\cal O}({\rm lg} N)$.

Related benchmarks for $f_{\rm F3d}$ by Bozzini, Mira and Rossini found  $E_{\tt rms} = 4\times 10^{-4}$ 
and $E_{\infty} = 8\times 10^{-3}$ using $N_p=3 375$ \cite{bozzini2002testing}, while Cavoretto, Rossi and Perracchione \cite{cavoretto2015} 
found $E_{\tt rms} = 6\times 10^{-5}$ and $E_{\infty} = 6\times 10^{-3}$ with $N_p=41 665$.  In this study,
it was possible to achieve compact support with fewer than $10^3$ octree nodes, $E_{\tt rms} < 10^{-11}$ and 
$E_{\infty} < 10^{-8}$; a roughly 5 order of magnitude improvement. 

\section{Conclusions}
In this contribution, a free-mesh interpolation was developed based on factorization of the Gaussian-RBF kernel in the 
flat polynomial limit, corresponding to Taylor expansion and the Vandermonde basis of geometric moments.
A top-down octree-scoping based on low-rank polynomial approximation was developed and demonstrated, achieving 
roughly 5 orders of magnitude improvement in free-mesh interpolation errors for the three-dimensional 
Franke function, relative to previous benchmarks.  A small advantage of the current implementation over
the continuous Gaussian RBF interpolation is the ease of implementation, and the potential for a simplified 
downward propagation of smooth components, {\em e.g.} in implementations with a scoping grid. 

This contribution also highlighted interesting forms of the identity, viz Eq.(\ref{identityprime}), 
obtaining also from the flat limit of the continuous case \cite{steinwart2006explicit,fasshauer2011positive}:
\begin{equation*}
\mat{I}_{\cal H} =  \left|\mat{\sqrt{\lambda}\varphi}' \right> \left< \mat{\varphi/\sqrt{\lambda}} \right| \, .
\end{equation*}
Interpolation in this case achieves super-spectral rates of convergence \cite{fornberg2005accuracy,fasshauer2011positive},
enabling many additional orders of precision to be obtained in the tuned shape-parameter ``dip'' \cite{fornberg2005accuracy,fasshauer2011positive}.  

Recently, we developed an $n$-body method for ill-conditioned Gramian factorization through octree-scoping on
the sub-space metric of the product manifold of the matrix-matrix multiply, achieving multiple 
order of magnitude compression through 
strongly contractive identity iteration \cite{Challacombe2015}. In that effort, a nested product form of the 
factor and its inverse was advocated; these factored forms are certainly compatible with low-rank factors alternatively found through 
QR factorization, via methods motivated by this and related work. Hybrid formulations of these two fast approaches 
are an interesting possibility, with the potential for extreme kernel factorizations. 

\bibliographystyle{siam}
\bibliography{FreeMesh}

\end{document}